\documentclass[11pt]{article}

\usepackage{amsfonts, amsmath, amssymb, amsthm}
\usepackage[a4paper,left=25mm,right=25mm,top=24mm,bottom=24mm]{geometry}
\usepackage{latexsym}
\usepackage{enumerate}
\usepackage{verbatim}
\usepackage{graphicx}
\usepackage[all,2cell,ps]{xy}
\usepackage{setspace} 
\usepackage{titlesec}
\usepackage[titles]{tocloft}
\usepackage{ragged2e}
\usepackage{indentfirst}
\usepackage[utf8]{inputenc}
\usepackage{tocloft}
\usepackage{xspace}

\usepackage{titlesec}
\usepackage[breaklinks]{hyperref}
\usepackage{blindtext}
 \usepackage{ifthen}
\usepackage{paralist}
\usepackage{makeidx}
\usepackage{nomentbl}
\usepackage{booktabs}
\usepackage[utf8]{inputenc}

\renewcommand{\theequation}{\arabic{section}.\arabic{equation}}

\newfam\msbmfam
\font\twelvemsbm=msbm10 scaled
\textfont\msbmfam\twelvemsbm\scriptfont\msbmfam\tenmsbm

\RaggedRight

\newtheoremstyle{newthm} 
     {3pt}
     {3pt}
     {}
     {}
     {\scshape}
     {:}
     {.5em}
     {}

\theoremstyle{newthm}

\usepackage[none]{hyphenat}
\usepackage{titlesec}
\usepackage{setspace} 
\titleformat*{\section}{\normalsize \bfseries}
\usepackage{ifpdf} 
 
\begin{document}


\begin{center}
    \textbf{Generalization Of Separation Of Variables  n-Harmonic Equation m Dimension and Unbounded Boundary Value Problem}
\end{center}

\begin{center}
    
    Ibraheem Otuf \\
    Department of mathematics \\
    Northern Border Univeristy  \\
    Saudi Arabia \\
    e-mail:Ibrahem.Otuf@nbu.edu.sa  \\
    
\end{center}
  



\singlespacing
\bigskip\bigskip

\noindent{\textbf{ABSTRACT}} The method of separation of variables is significant, it has been applied to physics, engineering , chemistry and other fields. It allows to reduce the diffculity of problems by separating the variables from partial differential equation system into ordinary differential equations system. However, this method has complexity in higher order partial differential equations. In this reserach, we generalize this method by using multinomial theorem of n-harmonic equation to solve n-harmonic equation with $m$ dimension and then solving an important class of partial differential equations with unbounded boundary conditions. Additionaly, application of convolution.     \\

  \def\theequation{\arabic{equation}}
        

        \section{Introduction}
        The n-harmonic 2 dimesion has been studied by Doscharis in \cite{Doscharis}, and  neumann Problem 1-harmonic m dimesion have been solved by Otuf in \cite{Otuf}. We are intersted in n-harmonic equation m dimension which is in the form of 
           \begin{equation}  
           \triangle_{m}^{n} u(x)=0,  
           \end{equation}
           where $  x=(x_1,x_2,........, x_m) . $
           This equations has high complexity since it is a partial differential equation of order 2n. However, separation of variables ia an intersting approch to solve this kind of problem. At begining,  we need to put our attention in driving new concept which decrease the complexity of solving this problem. Multinomial theorem of n-harmonic equation which is coming directly from the original multinomial theorem allows to wrtie the n-harmonic equation as summuation of ordinary differential equations. This lead to separate these equations and write them as sequnce of equations. Solving these equations will desgin an important gereralization. Moreover, we are intersted in solving unbounded boundary value problem that satisfying $ \triangle_m^{n} u(x)=0 $ and then application of its convolutions theorem.

           \section{ Multinomial Theorem of n-Harmonic Equation}
           
          We formalize multinomial theorem to desgin short approch to separate the variables of $ \triangle_m^{n} u(x) ,$  and 
          write them as summutions of ordinary differential equations. Supose for   $ 
          y_i=  \frac{X_{i}^{(2)}}{X_{i}} ,$  we have $ y_i^{(n)} = \frac{X_{i}^{(2n)}}{  X_i } $ . We let
 $  u(x)= \prod_{i = 1}^{m-1} X_i (x_i)   X_m (x_m) ,$ then apply $ \triangle_m^{n}  $ to  $ u(x) $ and then divide by $ \prod_{i = 1}^{m-1} X_i(x_i)   X_m (x_m)  $  leads to the  multinomial theorem of n-harmonic equation which is as following
          
          \begin{equation}
           \left(   \sum_{r=1}^{m} y_i \right)^{(n)}  = \sum_{h_1 + h_2+......+h_m=n}  \binom {n} {h_1, h_2, ......h_m}  \prod_{i = 1}^{m} y_i^{h_i}.
          \end{equation} 
          Where  
          \begin{equation}
          \binom {n} {h_1, h_2, ...... , h_m}  =   \frac{n!}{h_1 ! h_2 ! .........h_m !}, 
          \end{equation} 
           $ (n) $ is order of dervative and $ h_i \neq 0. $  
          Furthermore, mutiplication of any two or more ordinary differential equations, yields  sum of their orders. For example 
          $ X_{2}^{(2)} X_{2}^{(4)}= X_{2}^{(6)} $  
      We list some examples 
           For $ \triangle_2^{2} u(x)=0 , $ we have
          \begin{equation} ( y_1 + y_2 )^{(2)} =   y_1^{(2)}+ 2 y_1y_2  +  y_2^{(2)} = \frac{X_{1}^{(4)}}{X_{1}}  
          + 2 \frac{X_{1}^{(2)}}{X_{1}} \frac{X_{2}^{(2)}}{X_{2}} + \frac{X_{2}^{(4)}}{X_{2}}   .
          \end{equation}
          For $ \triangle_3^{2} u(x)=0 , $ we have 
     \begin{equation}      
     ( y_1 + y_2  )^{(3)}= y_1^{(3)}+  3 y_1^{(2)} y_2  +  3 y_1 y_2^{(2)}  + y_2^{(3)} =\frac{X_{1}^{(6)}}{X_{1}} + 3 \frac{X_{1}^{(4)}}{X_{1}}  \frac{X_{2}^{(2)}}{X_{2}} + 3 \frac{X_{1}^{(2)}}{X_{1}} \frac{X_{2}^{(4)}}{X_{2}} + \frac{X_{1}^{(6)}}{X_{1}}  .
       \end{equation}
          For $ \triangle_3^{3} u(x)=0 , $ we have 
       \begin{equation}      
    ( y_1 + y_2  + y_3)^{(3)}=   \frac{X_{1}^{(6)}}{X_{1}} + \frac{X_{2}^{(6)}}{X_{2}} +  \frac{X_{3}^{(6)}}{X_{3}}
          +  3 \frac{X_{1}^{(4)}}{X_{1}}  \frac{X_{2}^{(2)}}{X_{2}} + 3 \frac{X_{1}^{(2)}}{X_{1}} \frac{X_{2}^{(4)}}{X_{2}}   \end{equation} 
           \begin{equation}   +  3 \frac{X_{1}^{(4)}}{X_{1}}  \frac{X_{3}^{(2)}}{X_{3}} + 3 \frac{X_{1}^{(2)}}{X_{1}} \frac{X_{3}^{(4)}}{X_{3}} +  3 \frac{X_{3}^{(4)}}{X_{3}}  \frac{X_{2}^{(2)}}{X_{2}} + 3 \frac{X_{3}^{(2)}}{X_{3}} \frac{X_{2}^{(4)}}{X_{2}}  .
          \end{equation}

\section{Generalization Of Separation Of Variables}


           Let    $ x=(x_1,x_2,........x_m) $ be vector in $ \mathbb{R}^{m}. $ The n-harmonic equation is
             \begin{equation}  
           \triangle_m^{n} u(x)=0.  
           \end{equation}
                Let $ \lambda_{i_j} $ to be an constants in $ \mathbb{R},  $
             where $ i=1,2,3,.......m ,$ and $j=1,2,......,n.$
             We seek solution of the form
              \begin{equation}
               u(x)=\prod_{i = 1}^{m-1} X_i(x_i)   X_m (x_m) .    \label{eq-110}  
              \end{equation}
              Apply method of sepration of varibles to get multinomial theorem of n-harmonic equation that is  \begin{equation}
           \left(   \sum_{r=1}^{m} y_i \right)^{(n)}  = \sum_{h_1 + h_2+......+h_m=n}  \binom {n} {h_1, h_2, ......h_m}  \prod_{i = 1}^{m} y_i^{h_i} .    
          \end{equation}  
          
          Where  
          \begin{equation}
            y_i=  \frac{X_{i}^{(2)}}{X_{i}} , 
            \end{equation}
          
                 and 
                 \begin{equation}
                  y_i^{(n)} = \frac{X_{i}^{(2n)}}{ {X_{i}}} .
                  \end{equation}

          Let \begin{equation}
           \lambda_{i_j}= \frac{X_{i}^{(2j)}}{X_{i}} ,
           \end{equation} 
                       we get ordinary differential equations those are 
               \begin{equation}
               X_{i}^{(2j)} - \lambda_{i_j} X_{i}=0,      \label{eq-111}  
              \end{equation}
               and 
              \begin{equation} 
              \sum_{k=0}^{(n)} \binom {n} {k}  X_{m}^{(2n-2k)}  K_m^{(2k)}  = 0 .  \label{eq-112}  
              \end{equation}
              Where   $ X_{m}^{(0)} = X ,      $
                $ (2j) $ is the order of dervative, and   $ k_m= \sum_{i =1}^{m-1} \lambda_{i_1}.  $  
                Let solve \eqref{eq-111} if  $ j=1, $    we get
                 \begin{equation} 
                 X_{i}^{(2)} - \lambda_{i_1} X_{i}=0.  \label{eq-113} 
                \end{equation}  
                 If $ \lambda_{i_1} < 0, $ we have
                 \begin{equation} 
                 X_i (x_i)= a_i \cos \sqrt{\lambda_{i_1}} x_{i}  + b_i \sin \sqrt{\lambda_{i_1}} x_{i}.   \label{eq-1}  
                  \end{equation}
                for  $ j > 1 . $  we will let $ \lambda_{i_j} = (\lambda_{i_1})^j   $ in   \eqref{eq-111}, then the solutions of 
                this forms have two terms. First one is  \eqref{eq-1}, and second terms in which their constant forced to be 
                zeros. The final result is \eqref{eq-1}.    
                Similary if $ \lambda_{i_1} > 0,$    we  have 
                 \begin{equation} 
                 X_i (x_i) = a_i \cosh \sqrt{\lambda_{i_1}} x_{i}  + b_i \sinh \sqrt{\lambda_{i_1}} x_{i}  .                   
                 \label{eq-2}  
                  \end{equation}
                  If $ \lambda_{i_1} = 0,$   we have 
                   \begin{equation} 
                   X_i(x_i) = a_i   + b_i x .    
                      \label{eq-3}  
                    \end{equation}
                  
                The solution of \eqref{eq-112}, elementary differetial equation yields
                 \begin{equation}  \sum_{k=0}^{n} ( X_{m}^{2} + k_m )^{(n)}=0 .
                 \end{equation} Where $ X^{(0)}=X . $  This gives  
                   \begin{equation} 
                   \sum_{k=0}^{n} ( s^{2} + k_m )^{n}=0 .
                   \end{equation}
                  If $ k_m < 0 ,$  we have 
                \begin{equation} 
                 X_m(x_m)= \sum_{r=1}^{2n} x_m^{r-1}[ c_r \cosh \sqrt{- k_m} x_m  +  d_r \sinh \sqrt{ - k_m} x_m ].  \label{eq-4}  
                \end{equation} 
                 If $ k_m > 0 ,$  we have 
                  \begin{equation}  
                  X_m(x_m)= \sum_{r=1}^{2n} x_m^{r-1}[ c_r \cos  \sqrt{ k_m} x_m)  +  d_r 
                 \sin  \sqrt{k_m}x_m ].  \label{eq-5}  
                \end{equation} 
                Finally if $ k_m = 0 ,$  we have 
                \begin{equation}
                  X_m(x_m)= \sum_{r=1}^{2m} c_r x^{r-1} .   \label{eq-6}  
                \end{equation} 
                 The final solutions depends on choice of  $ \lambda_{i_1} $ for each $ i=1,2,.........m-1. $  We have 
                   \begin{equation}
                    u(x)=\prod_{i = 1}^{m-1} X_i(x_i)   X_m (x_m) .    \label{eq-7}  
                    \end{equation}   
                Remark: We may have different signs for $ \lambda_{i_1} . $  
                \section{Boundary Value Problem Unbounded Domain}
    We consider the Boundry value problem (BVP) replace $ x=(x_1, x_2, ........, x_{m-1},x_m) $ by  $ (x,x_m)=(x_1, x_2, ........, x_{m-1},x_m).  $

    \begin{equation}
    \triangle_m^{n} u(x,x_m) =0 , \quad   - \infty < x_i <  \infty;  \quad  L < x_m < \infty , 
    \end{equation} 

 \begin{equation}
  u (x, 0)=f(x) ,  \  \  \mid u(x,x_m)\mid  \textless \  M  . 
   \end{equation}  
Write  $ u(x,x_m) =  \prod_{i = 1}^{m-1} X_i(x_i)   X_m (x_m). $
Apply  gererlazation of sepration of varibles to get 
          
    
                \begin{equation}
               X_{i}^{(2j)} - \lambda_{i_j} X_{i} =0 ,      
              \end{equation}
              and 
             \begin{equation}  
             \sum_{k=0}^{n} ( X_{m}^{2} + k_m )^{(n)}=0 .
                 \end{equation} 
                 Where  $ X^{(0)}=X .  $
                 From the gereraliztion we need only $ j=1 $ i.e. 
                 \begin{equation}
               X_{i}^{(2)} - \lambda_{i_1} X_{i}=0,      
              \end{equation} 
                      and by letting  $ \lambda_{i_1}= - w_{i}^2 .$ 
                     We get   \begin{equation}
               X_{i}^{(2)} + w_{i}^{2} X_{i}=0 .     
              \end{equation} 
              The solutions is 
        \begin{equation}
                      X_i (x_i)= a_i \cos w_i x_{i}  + b_i \sin w_i x_{i} .   
                    \end{equation} 
                    The solutions for the mth variable is in the form 
                    
                    \begin{equation}  
                    X_m(x_m)= \sum_{r=1}^{2n} x_m^{r-1}[ q_r e^{ - \sqrt{ - k_m} x_m}  +  f_r  e^{ \sqrt{ - k_m} x_m}  ].   
                    \end{equation} 
                    
                    The final solutions 
                  \begin{equation} 
                   u(x,x_m) =  \prod_{i = 1}^{m-1} [a_i \cos w_i x_{i}  + b_i \sin w_i x_{i} ]   \sum_{r=1}^{2n} x_m^{r-1}[ q_r e^{ - \sqrt{ - k_m} x_m}  +  f_r  e^{ \sqrt{- k_m} x_m}  ] .
                   \end{equation}  
                 Where $ k_m= - \sum_{i =1}^{m-1} w_i^{2} . $  Since  $ w_i  > 0  , $   and  $ x_m \rightarrow  \ \infty ,   $  we must have $ f_r = 0 .$  It follows that 
              \begin{equation} 
                u(x,x_m) =  \int_{(0,\infty)^{m-1} }  \sum_{r=1}^{2n}  \prod_{i = 1}^{m-1} [A_{i_r} x_m^{r-1} \cos w_i x_{i}  
              +  x_m^{r-1}  B_{i_r} \sin w_i x_{i} ] [ e^{ - \sqrt{ - k_m} x_m}   ] \   dw_1 dw_2 \cdots dw_{m-1} .
           \end{equation}  
                    
                    Since    $ u (x, L)=f(x) ,  $  then
                           \begin{equation} 
                            f(x)  =    \int_{(0,\infty)^{m-1} }\sum_{r=1}^{2n}   \prod_{i = 1}^{m-1} L^{r-1} e^{ - \sqrt{- k_m} L }    [ A_{i_r}   \cos w_i x_{i}  + A_{i_r}  \sin w_i x_{i} ]   \  d w_1 dw_2 \cdots  d w_{m-1} . 
                            \end{equation}  
                              Where $ w=(w_1, w_2, ........, w_{m-1}).  $
                    Thus, we have
                   \begin{equation}
                   A_{i_r} (w) =   \frac{ \sum_{r=1}^{2n} L^{1-r} e^{ \sqrt{ - k_m} L }}{  \pi^{m-1}} \int_{(-\infty,\infty)^{m-1} }  f(x)     \prod_{i = 1}^{m-1}  \cos w_i x_i \  d x_1 d x_2  \cdots   d x_{m-1} .
                 \end{equation} 
                     \begin{equation}
                     B_{i_r} (w) = \frac{\sum_{r=1}^{2n}  L^{1-r} e^{ \sqrt{- k_m} L } }{ \pi^{m-1}}  \int_{(-\infty,\infty)^{m-1} }  f(x)     \prod_{i = 1}^{m-1}  \sin w_i x_i \   d x_1 d x_2 \cdots      d x_{m-1} . 
                     \end{equation}   
                   Where      $ r=1,2,...........,2n .$

\newcommand{\N}{\mathbb{N}}
 
            \section{Application of Convolution Theorem}
             The case 1-harmonic 2 dimenstion has been done in \cite{Haberman} , we extended to  n-harmonic 2 dimenstion. 
             Let us the result of previous problem when m=2, 
              \begin{equation} 
               u(x_1, x_2) =  \int_0^{\infty}  \sum_{r=1}^{2n}   [A_{1_r} x_2^{r-1} 
              \cos w_1  x_{1}  +  x_2^{r-1}  B_{i_r} \sin w_1 x_{1} ] [ e^{ - w_1  x_2}   ] \ dw_1 .
               \end{equation}  
            where
           \begin{equation}   A_{1_r} (w) =   \frac{ \sum_{r=1}^{2n} L^{1-r} e^{ w_1 L }}{  \pi } \int_{-\infty}^{\infty}  f(x_1)      \cos w_1 x_1  \ dx_1 . 
            \end{equation} 
               \begin{equation}
                B_{1_r} (w) = \frac{\sum_{r=1}^{2n}  L^{1-r} e^{ w_1  L } }{ \pi } \int_{-\infty}^{\infty}   f(x_1)     \sin   w_1 x_1 \  dx_1.
                \end{equation} 
                After we plug the coefficients into the general formula, we get
                     \begin{equation} 
                      u(x_1, x_2) =  \frac{\sum_{r=1}^{2n}  L^{1-r} e^{ w_1} L }{ \pi}  \int_{-\infty}^{\infty}  \int_{0}^{\infty}   \sum_{r=1}^{2n}   [  x_2^{r-1} 
              \cos w_1  x_1 \cos w_1 z_1   f(z)  .   
                 \end{equation} 
              \begin{equation}   +   x_2^{r-1}   \sin w_1 x_{1} \sin w_1 z_1 f(z_1)   ] [ e^{ - w_1 x_2}   ] \ dw_1  dz_1 
              \end{equation} 
                 \begin{equation}  =  \frac{\sum_{r=1}^{2n}  L^{1-r} e^{ w_1 L } }{ \pi}  \int_{-\infty}^{\infty} \int_0^{\infty}   \sum_{r=1}^{2n}   [ x_2^{r-1} e^{ - w_1 x_2}  
               \cos w_1 (z_1 -x_1)   f(z_1)  ]  \ dw_1  dz_1 .
                \end{equation}  
               
               Let \begin{equation}
                g(x_1, x_2)=  \int_0^{\infty}  \sum_{r=1}^{2n}   \cos w_1 (x_1) x_2^{r-1} e^{ -  w_1 x_2}  \ dw_1  
                \end{equation}  
             \begin{equation}  =  \sum_{r=1}^{2n}    \int_0^{\infty}  \frac{e^{i w_1 x_1} +e^{-i w_1 x_1 }}{2}  e^{ - w_1 x_2} x_2^{r-1}   \ dw_1 
             \end{equation}  
           \begin{equation}    =  \sum_{r=1}^{2n}   x_2^{r-1} e^{ - w_1 x_2}  \frac{1}{2} \left( \frac{1}{x_2-ix_1} + \frac{1}{x_2+i x_1} \right)  .
            \end{equation}  
                \begin{equation}  
                        = \sum_{r=1}^{2n}   x_2^{r-1} e^{ - w_1 x_2}  \frac{x_2}{x_2^{2}+x_1^{2}} .          
                \end{equation}
           Define the convolution theorem in the form of 
               \begin{equation} 
                u(x_1,x_2) =  \frac{\sum_{r=1}^{2n}  L^{1-r} e^{ w_1 L } }{ \pi }  \int_{-\infty}^{\infty}  f(z_1)  g(x_1-z_1 , x_2 ) dz_1. 
                \end{equation}  
           Let          
   \begin{equation} 
             f(z_1)= \left \{ \begin{array}{rcl} 1  & \mbox{if} & z \textgreater  \ 0   , \\  0   &       \mbox{if} &  z \textless \ 0   . \\ \end{array}  \right \} 
    \end{equation} 
 
           Now apply the theorem to get 
           
           \begin{equation}  
           u(x_1,x_2) =  \frac{\sum_{r=1}^{2n}  L^{1-r} e^{ w_1 L } }{ \pi}   \int_{-\infty}^{\infty}   \sum_{r=1}^{2n}   x_2^{r-1} e^{ - w_1 x_2}  \frac{x_2}{x_2^{2}+(x_1 - z_1)^{2}}   dz_1. 
          \end{equation}  
             
           \begin{equation}  
            =  \frac{\sum_{r=1}^{2n}  L^{1-r} e^{ w_1 L } }{ \pi}    \sum_{r=1}^{2n}  x_2^{r-1} e^{ - w_1 x_m}  \int_{-\infty}^{\infty}   e^{ - \sqrt{ k_m} x_m}  \frac{x_m}{x_m^{2}+(x_i - z_i)^{2}}   dz_i. 
              \end{equation}  
                
                   \begin{equation}  
                    =  \frac{\sum_{r=1}^{2n}  L^{1-r} e^{ w_1 L } }{ \pi }    \sum_{r=1}^{2n} 
              x_2^{r-1} e^{ - w_1 x_2}  \left [  \frac{\pi }{ 2} + \tan^{-1} \left( 
             \frac{x_1}{x_2} \right) \right ]  .
              \end{equation}

             \section{Parabolic Version Of n-Harmonic Equation}
             The parabolic version of n-harmonic equation is  
           \begin{equation}
 \alpha \triangle_m^{n} u(x, t) = u_t (x, t).
   \end{equation} 

 Where $ (x,t)=(x_1,x_2,........x_m, t) , $ $ t >0 ,$  $ \alpha \in \mathbb{R} , $ and  $ (x,t) \in \mathbb{R}^m \times (0, \infty ) . $  
            The between  difference between the eillptic version of n harmonic equation and prarbolic one  is similar to difference between heat equation and  1-harmonic equation. The solution can presented as $  u(x)=\prod_{i = 1}^{m} X_i(x_i) T(t) .$   $X_i(x_i)  $ has same solutions as the one in n harmonic equation.   For $ T(t) $ we have 
           $$ T^{'} (t)- \alpha k_{m+1}T(t)=0 .$$   
           The solutions is $$ T(t)= A e^{ \alpha t k_{m+1}  } . $$
        
             \section{Hyperbolic Version Of n-Harmonic Equation}
             The hyperbolic version of n-harmonic equation is 
       \begin{equation}
   \beta^2 \triangle_m^{n} u(x, t) = u_t (x, t).
 \end{equation}  Where $ (x,t)=(x_1,x_2,........x_m, t) , $ $ t >0 ,$  $ \alpha \in \mathbb{R} , $ and  $ (x,t) \in \mathbb{R}^m \times (0, \infty ) . $  
             The  difference between the eillptic version of n harmonic equation and hyperbolic one is similar to difference between wave equation and  1-harmonic equation. The solution can be written  as $  u(x)=\prod_{i = 1}^{m} X_i(x_i) T(t) .$   $X_i(x_i)  $ has same solutions as the one in n harmonic equation.   For $ T(t) $ we have 
           $$ T^{''} (t)- \beta^2 k_{m+1} T(t)=0. $$   
           The solutions is  $$ T(t)= C e^{ \beta \sqrt {k_{m+1}} t }  +  D  e^{ - \beta t \sqrt {k_{m+1} t }}   . $$
        
        \section{Conclusion}
        
         Sepations of variables is powerful method that assist to reduce the diffcuilty of partial differetion equation problems and solve them. The n-harmonic m dimension has been generlaized. Moreover, unbounded boundary value problem has been solved and applications of its convolution. Hyperbolic version and parabolic version are dicussed.

\addcontentsline{toc}{section}{REFERENCES}

\singlespacing


 \end{document}